\documentclass[preprints,article,accept,moreauthors,10pt,a4paper]{my-mdpi}

\firstpage{1} 

\history{Submitted: Nov 4, 2018; Revised: Nov 25, 2018; Accepted to Axioms (ISSN 2075-1680): Nov 28, 2018}

\Title{Regional Enlarged Observability of Fractional Differential 
Equations with Riemann--Liouville\\ Time Derivatives}

\Author{Hayat Zouiten $^{1}$, Ali Boutoulout $^{2}$ and Delfim F. M. Torres $^{3,}$*\orcidC{}}

\AuthorNames{Hayat Zouiten, Ali Boutoulout and Delfim F. M. Torres}

\address{%
$^{1}$ \quad TSI Team, MACS Laboratory, Department of Mathematics and Computer Science,
Moulay Ismail University, Faculty of Sciences, 11201 Meknes, Morocco; zouiten.hayat1991@gmail.com\\
$^{2}$ \quad TSI Team, MACS Laboratory, Department of Mathematics and Computer Science, 
Moulay Ismail University, Faculty of Sciences, 11201 Meknes, Morocco; boutouloutali@yahoo.fr\\
$^{3}$ \quad Center for Research \& Development in Mathematics and Applications (CIDMA),\newline 
Department of Mathematics, University of Aveiro, 3810--193 Aveiro, Portugal; delfim@ua.pt}

\corres{Correspondence: delfim@ua.pt; Tel.: +351-234-370-668}

\abstract{We introduce the concept of regional enlarged observability 
for fractional evolution differential equations involving 
Riemann--Liouville derivatives. The Hilbert Uniqueness Method (HUM) 
is used to reconstruct the initial state between two prescribed functions, 
in an interested subregion of the whole domain, without the knowledge of the state.}



\begin{document}

\section{Introduction} 

The exact birthday of fractional calculus, and the idea of non-integer differentiation, 
goes back to the 17th century, precisely to September 30, 1695, 
when L'H\^opital wrote a question to Leibniz about the meaning of 
$\displaystyle \frac{\partial^{n}}{\partial t^{n}}$ 
in the case $n = \displaystyle \frac{1}{2}$ \cite[pp.~301--302]{LE}.
Since then, many great mathematicians have investigated around this question. 
We can mention the studies of Euler (1730), Lagrange (1772), Laplace (1812), Fourier (1822), Abel (1823), 
Liouville (1832), Riemann (1847), Greer (1859), Gr\"{u}nwald (1867), Laurent (1884), 
Heaviside (1892), Pincherle (1902), Marchaud (1927), Love (1938), Widder (1941), Riesz (1949), 
Feller (1952), among others \cite{SKM}. These mathematicians began to consider 
how to define a fractional derivative. In 1860s, Riemann and Liouville obtained 
now celebrated definitions of fractional operators, by extending 
the Cauchy integral formula. Such fractional operators 
have a major role in practical problems \cite{MR3657873}.
In particular, Heymans and Podlubny have shown that it is possible 
to attribute a physical meaning to initial conditions expressed in terms 
of Riemann--Liouville fractional derivatives on the field of viscoelasticity, 
which is more appropriate than standard initial conditions \cite{HP}.
From a mathematical point of view, fractional calculus is a generalization 
of the traditional differential calculus to integrals and derivatives of non-integer order. 

The fact that real systems are better described with non-integer order differential equations 
has attracted engineer's interest, making fractional calculus a tool used nowadays in almost 
every area of sciences. Indeed, in the last decades, fractional calculus has been recognized as one 
of the best tools to describe long-memory processes and materials, anomalous diffusion, 
long-range interactions, long-term behaviours, power laws, and allometric scaling \cite{MR3287328}. 
Such models are those described by differential equations containing fractional order derivatives. 
On the other hand, fractional calculus is especially efficient for modelling systems related 
to diffusion processes \cite{MR3571004}. In \cite{DD,SDSP}, the heat transfer process 
was successfully modelled using a fractional model based on normal and anomalous diffusion 
equations. In \cite{DDG,DS,DSS}, very accurate models for ultra-capacitors
and electrical energy storage elements based on diffusion and the 
Helmholtz effect are presented. Simultaneously, fractional calculus has played 
a very important role in various fields such as physics, chemistry, mechanics, 
electricity, economics, signal and image processing, biophysics, and bioengineering 
\cite{HI,MA,MR3714436}. Likewise, for control theory, fractional calculus has 
an enormous role \cite{MR3316531,MR2721980,XUE}. The reader interested in applications 
of fractional calculus in control and mathematical modelling of systems and processes in physics, 
aerodynamics, electrodynamics of complex medium, viscoelasticity, heat conduction, 
and electricity mechanics, is referred to \cite{AX,BG,KST,LH,LS,OU,PO,SKM} and references therein.
For numerical methods to fractional partial differential equations, 
see \cite{MR3381791,Li:Yi:Kurths,MR3783171}.

Despite its development, the theory of fractional differential equations, 
compared with the classical theory of differential equations, 
is a field of research only on its initial stage of development, 
calling great interest to many mathematicians \cite{MR2962045,MyID:387,MR3763885}.
For distributed parameter systems, several works deal with the problem of regional 
observability, which we study here in the fractional context, investigating
the possibility to reconstruct the initial state or gradient only on a subregion 
$\omega$ of the evolution domain $\Omega$ \cite{AM,DR,EL,JAI,ZB}. 
For results on controllability, see \cite{CK,MR3172421,MR3244467,MyID:391}.
The interest to study the concept of observability for fractional differential 
equations is not new: see \cite{GCK,SP,XUE,ZA}. 
Here we investigate, for the first time in the literature, 
the concept of regional enlarged observability, that is, 
observability with constraints on the state for fractional diffusion equations.
For that, we make use of the Hilbert Uniqueness Method (HUM) of Lions \cite{LI,JL}. 

The paper is organized as follows. In Section~\ref{sec:2}, we present 
the problem of regional enlarged observability of fractional diffusion 
systems with the traditional first-order time derivative replaced by 
the Riemann--Liouville time fractional derivative. In Section~\ref{sec:3}, 
we give some preliminary results, which will be used throughout the paper.
In Section~\ref{sec:4}, we characterize the enlarged observability of the system. 
Section~\ref{sec:5} is focused on the regional reconstruction of the initial 
state in an internal subregion of the evolution domain. We present
in Section~\ref{sec:6} an example to demonstrate our main results.
We end with Section~\ref{sec:7} of conclusions 
and some possible directions of future research.


\section{Problem statement}
\label{sec:2}

In this section we formulate the concept of regional enlarged observability 
for a Riemann--Liouville time fractional diffusion system of order $\alpha \in (0,1]$. 
Let $\Omega$ be an open bounded subset of  $\mathbb{R}^{n}$ ($n = 1,2,3$), 
with a regular boundary $\partial \Omega$. For $T > 0$, let us denote 
$Q_{T} = \Omega \times [0,T]$ and $\Sigma_{T} = \partial \Omega \times [0,T]$. 
We consider the following time fractional order diffusion system: 
\begin{equation}
\label{se11}
\begin{cases}
{}_{0}D_{t}^{\alpha}y(x,t) =  A y(x,t) & \hbox{in}  \quad  Q_{T}, \\
y(\xi,t) = 0 &  \hbox{on} \quad \Sigma_{T},\\ 
\lim\limits_{t \rightarrow 0^{+}} {}_{0}I_{t}^{1-\alpha} y(x,t) 
= y_{0}(x) &  \hbox{in} \quad \Omega,
\end{cases}
\end{equation}
where ${}_{0}D_{t}^{\alpha}$ and ${}_{0}I_{t}^{\alpha}$ denote, 
respectively, the (left) Riemann--Liouville fractional derivative 
and integral with respect to time $t$, with $\alpha \in \mathbb{R}$ 
such that $0< \alpha\leq 1$. For details on these operators, see, 
e.g., \cite{MR3443073,KST,PO,MR3381791}. Here we just recall their definition: 
$$
{}_{0}D_{t}^{\alpha}y(x,t) 
= \displaystyle \frac{\partial}{\partial t}\,  
{}_{0}I_{t}^{1-\alpha} y(x,t)
$$
and 
$$
{}_{0}I_{t}^{\alpha} y(x,t) 
= \displaystyle \frac{1}{\Gamma(\alpha)} \int_{0}^{t} (t-s)^{\alpha-1}y(x,s) ds,
$$
where $\Gamma (\alpha) $ denotes Euler's Gamma function.
The second order operator $A$ in (\ref{se11}) is linear with dense domain, 
such that the coefficients do not depend on time $t$ and generates a strongly 
continuous semi-group $(S(t))_{t\geq 0}$ on the Hilbert space $L^{2}(\Omega)$. 
We assume that the initial state $y_{0}\in L^{2}(\Omega)$ is unknown. 
The observation space is $\mathcal{O} = L^{2}(0,T;\mathbb{R}^{q})$.

Without loss of generality, we denote $y(\cdot, t) := y(t)$. 
The measurements are obtained by the output function given by
\begin{equation}
\label{se12}
z(t) = C y(t), \quad t \in [0,T] \, ,
\end{equation}
where $ C $ is called the observation operator, which is a linear operator 
(possibly unbounded),  depending on the structure and the number $q\in \mathbb{N}$ 
of the considered sensors, with dense domain 
$D(C) \subseteq L^{2}(\Omega)$ and range in $\mathcal{O}$.


\section{Preliminaries}
\label{sec:3}

In this section, we recall some results for Riemann--Liouville time fractional 
differential systems and some notions and results to be used thereafter.

\begin{Lemma}[See \cite{CK,LL,MPG}]
\label{lemma3.1}
For any $u_{0} \in L^{2}(\Omega)$, $0<\alpha \leq 1$, we say that 
the function $u \in L^{2}(0,T;L^{2}(\Omega))$ is a mild solution of the system
$$
\begin{cases}
{}_{0}D_{t}^{\alpha}u(t)=   A u(t),  &  t \in [0,T],\\
\lim\limits_{t \rightarrow 0^{+}} {}_{0}I_{t}^{1-\alpha} u(t) 
=  u_{0}, 
\end{cases}
$$
if $u$ satisfies the equation
$$
u(t) =  H_{\alpha}(t)u_{0},
$$	
where
$$
H_{\alpha}(t) = \alpha  t^{\alpha-1} \int_{0}^{\infty} 
\theta \xi_{\alpha} (\theta) S(t^{\alpha}\theta) d\theta,
$$
$$
\xi_{\alpha} (\theta) = \displaystyle \frac{1}{\alpha} 
\theta^{-1-\frac{1}{\alpha}}\varpi_{\alpha}(\theta^{-\frac{1}{\alpha}}),
$$
$$
\varpi_{\alpha}(\theta)  = \displaystyle \frac{1}{\pi} 
\sum_{n=1}^{\infty} (-1)^{n-1} \theta^{-n\alpha-1} 
\frac{\Gamma(n\alpha+1)}{n!} \sin(n\pi\alpha), 
\quad \theta \in(0,\infty),
$$
with $\xi_{\alpha}$ the probability density function defined on 
$(0,\infty)$, satisfying 
$$
\xi_{\alpha}(\theta) \geq 0, \quad \theta \in (0,\infty)  
\quad \hbox{and} \quad \int_{0}^{\infty} \xi_{\alpha}(\theta) d\theta = 1.
$$
Moreover,
$$
\int_{0}^{\infty} \theta^{\nu} \xi_{\alpha}(\theta) d\theta 
= \frac{\Gamma(1+\nu)}{\Gamma(1+\alpha \nu)}, \quad \nu \geq 0.
$$
\end{Lemma}

Note that a mild solution of system (\ref{se11}) can be written as
\begin{equation*}
y(x,t) = H_{\alpha}(t) y_{0}, \quad t \in [0,T].
\end{equation*}
In order to prove our results, the following lemma is used.

\begin{Lemma}[See \cite{MAG}]
Let the reflection operator $\mathcal{Q}$ on the interval $[0, T]$ be defined by
$$
\mathcal{Q} f(t) := f(T-t),
$$
for some function $f$ that is differentiable and integrable 
in the Riemann--Liouville sense. Then the following relations hold:
$$
\mathcal{Q} {}_{0}D_{t}^{\alpha} f(t)
= {}_{t}D_{T}^{\alpha} \mathcal{Q} f(t), 
\qquad \mathcal{Q} 	{}_{0}I_{t}^{\alpha}   
f(t) = {}_{t}I_{T}^{\alpha} \mathcal{Q} f(t)
$$
and 
$$
{}_{0}D_{t}^{\alpha}\mathcal{Q} f(t) 
= \mathcal{Q} {}_{t}D_{T}^{\alpha} f(t), 
\qquad {}_{0}I_{t}^{\alpha} \mathcal{Q} f(t)
= \mathcal{Q}   {}_{t}I_{T}^{\alpha}   f(t).
$$
\end{Lemma}

Follows some notions of admissibility of the output operator $C$. 
The output function of the autonomous system (\ref{se11}) 
is expressed by
\begin{equation*}
z(t) = CH_{\alpha}(t)y_{0}  = K_{\alpha}(t)y_{0}, \quad t \in [0,T],
\end{equation*} 
where $K_{\alpha} : L^{2}(\Omega) \longrightarrow \mathcal{O}$ 
is a linear operator. To obtain the adjoint operator 
of $K_{\alpha}$, we have two cases.
\begin{description}
\item[Case 1.] $C$ is bounded (e.g., zone sensors).
Let $C : L^{2}(\Omega) \longrightarrow \mathcal{O}$ and 
$C^{*}$ be its adjoint. We get that the adjoint operator 
of $K_{\alpha}$ can be given by
$$ 
\begin{array}{rll}
K_{\alpha}^{*} : \mathcal{O} & \longrightarrow &  L^{2}(\Omega)\\
z^{*} & \longrightarrow& \displaystyle \int_{0}^{T} H_{\alpha}^{*}(s)C^{*}z^{*}(s)ds.
\end{array}
$$

\item[Case 2.] $C$ is unbounded (e.g., pointwise sensors).
In this case, we have 
$$
C : D(C) \subseteq L^{2}(\Omega) \longrightarrow \mathcal{O}
$$ 
with $C^{*}$ denoting its adjoint. In order to give a sense to (\ref{se12}), 
we make the assumption that $C$ is an admissible observation operator 
in the sense of Definition~\ref{def3}.
\end{description}

\begin{Definition}
\label{def3}
The operator $C$ of system (\ref{se11})--(\ref{se12}) is 
an admissible observation operator if there exists 
a constant $M > 0$ such that
\begin{equation*}
\displaystyle \int_{0}^{T} \left\| CH_{\alpha}(s)y_{0} \right\|^{2} ds 
\leq M \left\| y_{0} \right\|^{2} 
\end{equation*}
for any $y_{0} \in D(C)$.
\end{Definition}

Note that the admissibility of $C$ guarantees that we can extend the mapping 
$$
y_{0} \longmapsto CH_{\alpha}(t)y_{0} = K_{\alpha}(t)y_{0}
$$ 
to a bounded linear operator from 
$L^{2}(\Omega)$ to $\mathcal{O}$. For more details, see, e.g., \cite{PW,SA,WE}. 
Then the adjoint of the operator $K_{\alpha}$ can be defined as
$$ 
\begin{array}{rll}
K_{\alpha}^{*} :  D(K_{\alpha}^{*}) \subseteq \mathcal{O} 
& \longrightarrow &  L^{2}(\Omega)\\
z^{*} & \longrightarrow& 
\displaystyle \int_{0}^{T} H_{\alpha}^{*}(s)C^{*}z^{*}(s)ds.
\end{array}
$$ 


\section{Enlarged observability and characterization}
\label{sec:4}

Let $\omega$ be a subregion of $\Omega$ with a positive Lebesgue measure. 
We define the restriction operator $\chi_{\omega}$ 
and its adjoint $\chi_{\omega}^{*}$  by
\begin{equation*}
\begin{array}{rll}
\chi_{\omega} : L^{2}(\Omega)& \longrightarrow &  L^{2}(\omega)\\
y & \longrightarrow& \chi_{_{\omega}}y = y_{|\omega} \\
\end{array}
\end{equation*}
and
$$
(\chi_{_{\omega}}^{*}y)(x) 
= 
\begin{cases}
y(x) & \text{ if } \ x \in \omega,\\
0 & \text{ if } \ x \in \Omega\backslash\omega.
\end{cases}
$$
Similarly to the discussions in \cite{CZW,DR,PW}, it follows 
that a necessary and sufficient condition for the regional 
exact observability of the system described by (\ref{se11}) 
and (\ref{se12}) in $\omega$ at time $t$ is that 
$Im(\chi_{_{\omega}}K_{\alpha}^{*}) = L^{2}(\omega)$.
                                                              
Let $\beta(\cdot)$ and $\gamma(\cdot)$ be two functions defined 
in $L^{2}(\omega)$ such that $\beta(\cdot) \leq \gamma(\cdot)$ 
a.e. in $\omega$. Throughout the paper, we set
\begin{equation*}
[\beta(\cdot), \gamma(\cdot)] = \left\{ y\in L^{2}(\omega) \, | \, 
\beta(\cdot) \leq y(\cdot) \leq \gamma(\cdot)
\quad \hbox{a.e.}\; \hbox{in} \; \omega \; \right\}.
\end{equation*}
We consider
\begin{equation*}
y_{0} = 
\begin{cases}
y_{0}^{1}  & \hbox{in } \; [\beta(\cdot),\gamma(\cdot)], \\
y_{0}^{2} & \hbox{in } \; L^{2}(\Omega) \backslash [\beta(\cdot),\gamma(\cdot)].    
\end{cases}
\end{equation*}
The study of regional enlarged observability for the Riemann--Liouville 
time fractional order diffusion system amounts to solving the following problem.

\begin{Problem}
Given the system (\ref{se11}) together with the output (\ref{se12}) in $\omega$ 
at time $t \in [0,T] $, is it possible to reconstruct $y_{0}^{1}$ between 
two prescribed functions $\beta(\cdot)$ and $\gamma(\cdot)$ in $\omega$?
\end{Problem}

Before proving our first result, we need two important definitions.

\begin{Definition}
The system (\ref{se11}) together with the output (\ref{se12}) is said to be exactly 
$[\beta(\cdot), \gamma(\cdot)]$--observable in $\omega$ if 
\begin{equation*}
Im (\chi_{_{\omega}} K_{\alpha}^{*}) 
\cap [\beta(\cdot), \gamma(\cdot)] \neq \emptyset.
\end{equation*}
\end{Definition}

\begin{Definition}
The sensor $(D,f)$ is said to be exactly $[\beta(\cdot), \gamma(\cdot)]$--strategic 
in $\omega$ if the observed system is exactly 
$[\beta(\cdot), \gamma(\cdot)]$--observable in $\omega$.
\end{Definition}

\begin{Remark}
If $\alpha = 1$, then system (\ref{se11}) is reduced to the
normal diffusion process recently considered in \cite{ZO}.
The results of \cite{ZO} are a particular case of our results.
\end{Remark}
	
\begin{Remark}
If the system (\ref{se11}) together with the output (\ref{se12}) is exactly 
$[\beta(\cdot), \gamma(\cdot)] $--observable in $\omega_{1}$, then it is 
exactly $[\beta(\cdot), \gamma(\cdot)]$--observable in any subregion 
$\omega_{2} \subset \omega_{1}$.
\end{Remark}

\begin{Theorem}
\label{thm4.1}
The following two statements are equivalent:
\begin{enumerate}
\item[\textbf{1.}] The system (\ref{se11}) together 
with the output (\ref{se12}) is exactly 
$[\beta(\cdot), \gamma(\cdot)]$--observable in $\omega$.

\item[\textbf{2.}] $Ker (K_{\alpha} \chi_{_{\omega}}^{*}) 
\cap [\beta(\cdot), \gamma(\cdot)]  = \{0\}$.
\end{enumerate}
\end{Theorem}

\begin{proof}
We begin by proving that statement 1 implies 2.
For that we show that
$$
\begin{array}{rll}
Im (\chi_{_{\omega}} K_{\alpha}^{*}) 
\cap [\beta(\cdot), \gamma(\cdot)] \neq \emptyset
&\Longrightarrow
& Ker (K_{\alpha} \chi_{_{\omega}}^{*}) 
\cap [\beta(\cdot), \gamma(\cdot)]  = \{0\}.
\end{array}
$$
Suppose that
$$
Ker (K_{\alpha} \chi_{_{\omega}}^{*}) 
\cap [\beta(\cdot), \gamma(\cdot)] \neq \{0\}.
$$
Let us consider $y \in Ker (K_{\alpha} \chi_{_{\omega}}^{*}) \cap[\beta(\cdot), \gamma(\cdot)]$ 
such that $y \neq 0$. Then, $y \in Ker (K_{\alpha} \chi_{_{\omega}}^{*}) $ and 
$y \in [\beta(\cdot), \gamma(\cdot)]$.
We have $Ker (K_{\alpha} \chi_{_{\omega}}^{*}) = Im(\chi_{_{\omega}} K_{\alpha}^{*})^{\perp}$, 
so that $y \in Im(\chi_{_{\omega}} K_{\alpha}^{*})^{\perp}$, $y \neq 0$. 
Therefore, $y \notin Im(\chi_{_{\omega}} K_{\alpha}^{*})$, and 		
\begin{gather*}
Ker (K_{\alpha} \chi_{_{\omega}}^{*}) \cap[\beta(\cdot), \gamma(\cdot)]  
\subset \displaystyle L^{2}(\omega)\setminus Im(\chi_{_{\omega}} K_{\alpha}^{*}),\\
Im(\chi_{_{\omega}} K_{\alpha}^{*})\subset \displaystyle 
\left[ L^{2}(\omega)\setminus Ker (K_{\alpha} \chi_{_{\omega}}^{*})\right]  
\cup \displaystyle \left[ L^{2}(\omega)\setminus [\beta(\cdot), \gamma(\cdot)] \,\right].
\end{gather*}
We have 
$$
Im(\chi_{_{\omega}} K_{\alpha}^{*}) \subset \displaystyle L^{2}(\omega) 
\setminus Ker (K_{\alpha} \chi_{_{\omega}}^{*}).
$$
Accordingly,
$$
Im(\chi_{_{\omega}} K_{\alpha}^{*}) 
\cap Ker (K_{\alpha} \chi_{_{\omega}}^{*}) = \emptyset
$$
and
$$
Im(\chi_{_{\omega}} K_{\alpha}^{*}) \cap Im(\chi_{_{\omega}} 
K_{\alpha}^{*})^{\perp} = \emptyset,
$$ 
which is absurd. Since
$$
Im(\chi_{_{\omega}} K_{\alpha}^{*}) \subset \displaystyle 
\displaystyle  L^{2}(\omega)\setminus [\beta(\cdot), \gamma(\cdot)],
$$
it follows that
$$
Im(\chi_{_{\omega}} K_{\alpha}^{*}) \cap [\beta(\cdot), \gamma(\cdot)]  
= \emptyset,
$$ 
which is also absurd. Consequently, 
$$
Ker (K_{\alpha} \chi_{_{\omega}}^{*}) \cap [\beta(\cdot), \gamma(\cdot)]  = \{0\}.
$$
We now prove the reverse implication: statement 2 implies 1. For that we show that
$$
\begin{array}{rll}
Ker (K_{\alpha} \chi_{_{\omega}}^{*}) 
\cap [\beta(\cdot), \gamma(\cdot)] = \{0\}
&\Longrightarrow& Im(\chi_{_{\omega}} K_{\alpha}^{*}) 
\cap [\beta(\cdot), \gamma(\cdot)] \neq \emptyset.
\end{array}
$$
Suppose that
$$
Ker (K_{\alpha} \chi_{_{\omega}}^{*}) 
\cap [\beta(\cdot), \gamma(\cdot)]  = \{0\}.
$$
Let us consider 
$$
y \in Ker (K_{\alpha} \chi_{_{\omega}}^{*})  \cap[\beta(\cdot), \gamma(\cdot)].
$$
Then, $y \in Ker (K_{\alpha} \chi_{_{\omega}}^{*})$ and 
$y \in [\beta(\cdot), \gamma(\cdot)]$ such that $y = 0$.
We have 
$$
Ker (K_{\alpha} \chi_{_{\omega}}^{*}) = Im(\chi_{_{\omega}} 
K_{\alpha}^{*})^{\perp}, 
$$
so $y \in Im(\chi_{_{\omega}} K_{\alpha}^{*})^{\perp}$ such that $y = 0$. Hence, 
$$
y \in Im(\chi_{_{\omega}} K_{\alpha}^{*})\; \hbox{and}\; y \in [\beta(\cdot), \gamma(\cdot)] 
$$ 
and
$$
Im(\chi_{_{\omega}} K_{\alpha}^{*}) \cap [\beta(\cdot), \gamma(\cdot)]  \neq \emptyset,
$$
which shows that (\ref{se11})--(\ref{se12}) is exactly 
$[\beta(\cdot), \gamma(\cdot)]$--observable in $\omega$.
\end{proof}


\section{The HUM approach}
\label{sec:5}

The purpose of this section is to present an approach that allows us to reconstruct 
the initial state of the system (\ref{se11}) between two prescribed functions 
$\beta(\cdot)$ and $\gamma(\cdot)$ in $\omega$. Our approach constitutes 
an extension of the Hilbert Uniqueness Method (HUM) developed by Lions \cite{LI,JL}. 
In what follows, $\mathcal{G}$ is defined by
\begin{equation}
\label{se19}
\mathcal{G} = \left\{g \in L^{2}(\Omega) \; | \; g = 0 \; 
\hbox{in} \; L^{2}(\Omega) \backslash [\beta(\cdot), \gamma(\cdot)]\right\}.
\end{equation}


\subsection{Pointwise sensors}

Let us consider system (\ref{se11}) observed by a pointwise sensor 
$(b,\delta_{b})$, where $b \in \overline{\Omega}$ is the sensor 
location and $\delta$ is the Dirac mass concentrated in $b$. 
For details on pointwise sensors we refer the reader to \cite{EL}. 
Here the output function is given by
\begin{equation}
\label{se110}
z(t) = \varphi(b,T-t), \quad t \in [0,T].
\end{equation}
For $\varphi_{0} \in \mathcal{G}$, we consider the following system:
\begin{equation}
\label{se111}
\begin{cases}
{}_{0}D_{t}^{\alpha}\varphi(x,t)
=  A \varphi(x,t) & \hbox{in}  \quad  Q_{T}, \\
\varphi(\xi,t) = 0 &  \hbox{on} \quad \Sigma_{T},\\ 
\lim\limits_{t \rightarrow 0^{+}} {}_{0}I_{t}^{1-\alpha} \varphi(x,t) 
= \varphi_{_{0}}(x) &  \hbox{in} \quad \Omega.
\end{cases}
\end{equation}
Without loss of generality, we denote $\varphi(x, t) := \varphi(t)$. 
System (\ref{se111}) admits a unique solution 
$\varphi \in L^{2}(0,T;H_{0}^{1}(\Omega)) \cap C(\Omega \times [0,T])$ 
given by $\varphi(t) =  H_{\alpha}(t)\varphi_{_{0}}$.
We consider a semi-norm on $\mathcal{G}$ defined by
\begin{equation}
\label{se112}
\varphi_{_{0}} \longmapsto \| \varphi_{_{0}} \|_{\mathcal{G}}^{2} 
= \displaystyle  \int_{0}^{T} \left\| C \varphi(T-t) \right\|^{2}  dt.
\end{equation}
The following result holds.

\begin{Lemma}
\label{lema5.1}
If the system (\ref{se11}) together with the output (\ref{se110}) 
is exactly $[\beta(\cdot), \gamma(\cdot)]$--observable in $\omega$, 
then (\ref{se112}) defines a norm on $\mathcal{G}$.
\end{Lemma}

\begin{proof}
Consider $\varphi_{_{0}} \in \mathcal{G}$. Then,
$$
\left\| \varphi_{_{0}} \right\|_\mathcal{G} 
\quad \Longrightarrow \quad C\varphi(T-t) = 0 
\quad \hbox{for all}  \quad t \in [0,T].
$$ 
We have 
$$
\varphi_{_{0}} \in L^{2}(\Omega) \quad \Longrightarrow 
\quad \chi_{_{\omega}}\varphi_{_{0}} \in L^{2}(\omega)
$$
or 
$$
K_{\alpha}(t) \chi_{_{\omega}}^{*}\chi_{_{\omega}}\varphi_{_{0}} 
= CH_{\alpha}(t) \chi_{_{\omega}}^{*}\chi_{_{\omega}}\varphi_{_{0}} = 0.
$$
Hence, 
$$
\chi_{_{\omega}}\varphi_{_{0}} \in Ker(K_{\alpha}\chi_{\omega}^{*}).
$$
For $\chi_{_{\omega}} \varphi_{_{0}} \in [\beta(\cdot), \gamma(\cdot)]$, 
one has $\chi_{_{\omega}} \varphi_{_{0}} \in Ker(K_{\alpha} \chi_{\omega}^{*}) 
\cap [\beta(\cdot), \gamma(\cdot)] $ and, because the system is exactly 
$[\beta(\cdot), \gamma(\cdot)]$--observable in $\omega$, $\chi_{_{\omega}}\varphi_{_{0}} = 0$. 
Consequently, $\varphi_{_{0}} = 0$ and (\ref{se112}) is a norm.
\end{proof}

Consider the system 
\begin{equation}
\label{se113}
\begin{cases}
\mathcal{Q} \; {}_{t}D_{T}^{\alpha}\Psi(x,t)
=  A^{*} \mathcal{Q} \Psi(x,t) + C^{*}C \mathcal{Q} \varphi(x,t)
& \hbox{in} \quad  Q_{T}, \\
\Psi(\xi,t) = 0  & \hbox{on} \quad \Sigma_{T},\\ 
\lim\limits_{t \rightarrow T^{-}} \mathcal{Q} \; {}_{t}I_{T}^{1-\alpha} \Psi(x,t) 
= 0 &  \hbox{in} \quad \Omega,
\end{cases}
\end{equation}
controlled by the solution of (\ref{se111}). For $\varphi_{_{0}} \in \mathcal{G}$, 
we define the operator $\Lambda : \mathcal{G} \longrightarrow \mathcal{G}^{*}$ by
\begin{equation*}
\Lambda \varphi_{_{0}} =  \mathcal{P}(\Psi(0)),
\end{equation*}  
where $\mathcal{P} = \chi_{\omega}^{*}\chi_{\omega}$ and $\Psi(0) = \Psi(x,0)$.
Let us now consider the system 
\begin{equation}
\label{se115}
\begin{cases}
\mathcal{Q} \; {}_{t}D_{T}^{\alpha}\Theta(x,t)
=  A^{*} \mathcal{Q}  \Theta(x,t) + C^{*} \mathcal{Q}z(t)
& \hbox{in} \quad  Q_{T}, \\
\Theta(\xi,t) = 0 &  \hbox{on} \quad \Sigma_{T},\\ 
\lim\limits_{t \rightarrow T^{-}}  
\mathcal{Q} \; {}_{t}I_{T}^{1-\alpha} \Theta(x,t) 
= 0 &  \hbox{in} \quad \Omega.
\end{cases}
\end{equation}
If $\varphi_{_{0}}$ is chosen such that $\Theta(0) = \Psi(0)$ in $\omega$, 
then system (\ref{se115}) can be seen as the adjoint of system 
(\ref{se11}) and our problem of enlarged observability is reduced to solve the equation
\begin{equation}
\label{se116}
\Lambda \varphi_{_{0}} = \mathcal{P}(\Theta(0)).
\end{equation} 

\begin{Theorem}
\label{thm5.1}
If system (\ref{se11}) together with the output (\ref{se110}) is exactly 
$[\beta(\cdot), \gamma(\cdot)] $--observable in $\omega$, then equation 
(\ref{se116}) admits a unique solution $\varphi_{_{0}} \in \mathcal{G}$, 
which coincides with the initial state $y_{0}^{1}$ observed between 
$\beta(\cdot)$ and $\gamma(\cdot)$ in $\omega$. Moreover,
$y_{0}^{1} = \chi_{_{\omega}} \varphi_{_{0}}$.
\end{Theorem}

\begin{proof}
By Lemma~\ref{lemma3.1}, we see that $\left\| \cdot \right\|_{\mathcal{G}}$ 
is a norm of the space $\mathcal{G}$ provided that the system (\ref{se11}) 
together with the output (\ref{se110}) is exactly 
$[\beta(\cdot), \gamma(\cdot)]$--observable in $\omega$.
Now, we show that (\ref{se116}) admits a unique solution in $\mathcal{G}$. 
For any $\varphi_{_{0}} \in \mathcal{G}$, equation (\ref{se116}) admits 
a unique solution if $\Lambda$ is an isomorphism. Then, 
$$
\begin{array}{rll}
\left\langle \Lambda \varphi_{_{0}}, \varphi_{_{0}} \right\rangle_{L^{2}(\Omega)}  
& = & \left\langle \mathcal{P}\Psi(0), \varphi_{_{0}} \right\rangle_{L^{2}(\Omega)} \\
& = & \left\langle \chi_{_{\omega}}^{*}\chi_{_{\omega}} 
\Psi(0), \varphi_{_{0}} \right\rangle_{L^{2}(\Omega)}  \\
& = & \left\langle \Psi(0), \varphi_{_{0}} \right\rangle_{L^{2}(\omega)}
\end{array}
$$
or $\Psi(t)$ is the solution of system (\ref{se113}), that is,
$$
\Psi(t) = H^{*}_{\alpha}(T-t)\Psi(T) 
+ \int_{t}^{T} H^{*}_{\alpha}(T-\tau)C^{*}C\varphi(T-\tau) d\tau
$$
and
$$
\Psi(0) = \displaystyle \int_{0}^{T} H^{*}_{\alpha}(T-\tau) 
C^{*}C\varphi(T-\tau) d\tau,
$$
where
$$
H^{*}_{\alpha}(t) = \alpha  t^{\alpha-1} \int_{0}^{\infty} 
\theta \xi_{\alpha} (\theta) S^{*}(t^{\alpha}\theta) d\theta
$$
with $(S^{*}(t))_{t \geq 0}$ the strongly continuous semi-group generated by $A^{*}$. 
We obtain that
$$
\begin{array}{rll}
\left\langle \Lambda \varphi_{_{0}}, \varphi_{_{0}} \right\rangle_{L^{2}(\Omega)}  
& = & \left\langle \Psi(0), \varphi_{_{0}} \right\rangle\\
& = & \displaystyle \left\langle  \int_{0}^{T} H^{*}_{\alpha}(T-\tau) 
C^{*}C\varphi(T-\tau) d\tau, \varphi_{_{0}} \right\rangle\\
& = & \displaystyle  \int_{0}^{T} \displaystyle \left\langle 
C\varphi(T-\tau), C H_{\alpha}(T-\tau) \varphi_{_{0}} \right\rangle d\tau\\
& = & \displaystyle \int_{0}^{T} \displaystyle \left\| 
C \varphi(T-\tau) \right\|^{2} d \tau\\
& = & \displaystyle \left\| \varphi_{_{0}} \right\|^{2}_{\mathcal{G}},
\end{array}
$$
concluding that $\Lambda$ is an isomorphism. Consequently, 
equation (\ref{se116}) has a unique solution that is also 
the initial state to be estimated between $\beta(\cdot)$ and 
$\gamma(\cdot)$ in the subregion $\omega$ given by
$$
y_{0}^{1} = \chi_{_{\omega}} \varphi_{_{0}}.
$$
The proof is complete.
\end{proof}


\subsection{Zone sensors}

Let us come back to system (\ref{se11}) and suppose that the measurements 
are given by an internal zone sensor defined by $(D,f)$ with $D \subset \Omega$ 
and $f \in L^{2}(D)$. The system is augmented with the output function
\begin{equation}
\label{se117}
z(t) = \displaystyle \int_{D} y(x,T-t) f(x) dx .
\end{equation}
In this case, we consider (\ref{se111}), $\mathcal{G}$ 
given by (\ref{se19}), and we define a semi-norm on $\mathcal{G}$ by
\begin{equation}
\label{se118}
\| \varphi_{_{0}}\|_{\mathcal{G}}^{2} 
= \displaystyle \int_{0}^{T} \displaystyle 
\left\langle  \varphi(T-t), f \right\rangle _{L^{2}(D)}^{2}dt
\end{equation}
with
\begin{equation*}
\begin{cases}
\mathcal{Q} \; {}_{t}D_{T}^{\alpha}\Psi(x,t)
=  A^{*} \mathcal{Q} \Psi(x,t) + \left\langle 
\mathcal{Q}\varphi(t), f \right\rangle_{L^{2}(D)} 
\chi_{_{D}}f(x)  & \hbox{in} \quad  Q_{T}, \\
\Psi(\xi,t) = 0 &  \hbox{on} \quad \Sigma_{T},\\ 
\lim\limits_{t \rightarrow T^{-}} \mathcal{Q} \; 
{}_{t}I_{T}^{1-\alpha} \Psi(x,t) = 0 
& \hbox{in} \quad \Omega.
\end{cases}
\end{equation*}
We introduce the operator
\begin{equation}
\label{se120}
\begin{array}{rll}
\Lambda : \mathcal{G}& \longrightarrow &  \mathcal{G}^{*}\\
\varphi_{_{0}} & \longrightarrow
& \Lambda \varphi_{_{0}} = \mathcal{P}(\Psi(0)),
\end{array}
\end{equation}   
where $\mathcal{P} = \chi_{\omega}^{*}\chi_{\omega}$ and $\Psi(0) = \Psi(x,0)$.
 Let us consider the system 
\begin{equation}
\label{se121}
\begin{cases}
\mathcal{Q} \; {}_{t}D_{T}^{\alpha}\Theta(x,t)
=  A^{*} \mathcal{Q} \Theta(x,t) + \left\langle \mathcal{Q} z(t), 
f \right\rangle _{L^{2}(D)} \chi_{_{D}}f(x) 
& \hbox{in} \quad  Q_{T}, \\
\Theta(\xi,t) = 0 &  \hbox{on} \quad \Sigma_{T},\\ 
\lim\limits_{t \rightarrow T^{-}} 
\mathcal{Q} \; {}_{t}I_{T}^{1-\alpha} \Theta(x,t) 
= 0 &  \hbox{in} \quad \Omega.
\end{cases}
\end{equation}
If $\varphi_{_{0}}$ is chosen such that $\Theta(0) = \Psi(0)$ in $\omega$, 
then (\ref{se121}) can be seen as the adjoint of system (\ref{se11}) 
and our problem of enlarged observability consists to solve the equation
\begin{equation}
\label{se122}
\Lambda \varphi_{_{0}} = \mathcal{P}(\Theta(0)).
\end{equation} 

\begin{Theorem}
If system (\ref{se11}) together with the output (\ref{se117}) 
is exactly $[\beta(\cdot), \gamma(\cdot)]$--observable in $\omega$, 
then equation (\ref{se122}) has a unique solution 
$\varphi_{_{0}} \in \mathcal{G}$, which coincides with the initial 
state $y_{0}^{1}$ observed between $\beta(\cdot)$ 
and $\gamma(\cdot)$ in $\omega$.
\end{Theorem}

\begin{proof}
The proof is similar to the proof of Theorem~\ref{thm5.1}.
\end{proof}


\section{Example}
\label{sec:6}

Let us consider the following one-dimensional time fractional
differential system of order $\alpha \in (0,1]$ in 
$\Omega_{1} = [0,1]$, excited by a pointwise sensor:
\begin{equation}
\label{se123}
\begin{cases}
{}_{0}D_{t}^{\alpha}y(x,t)
= \displaystyle \frac{\partial^{2}y(x,t)} {\partial x^{2}} 
& \hbox{in} \quad  [0,1]\times[0,T], \\
y(0,t) \;= \;y(1,t) = 0 &  \hbox{in} \quad [0,T],\\ 
\lim\limits_{t \rightarrow 0^{+}} {}_{0}I_{t}^{1-\alpha} y(x,t) 
= y_{0}(x) &  \hbox{in} \quad [0,1],
\end{cases}
\end{equation}
augmented with the output function
\begin{equation}
\label{se124}
z(t) = C y(x,t) = y(b,t),
\end{equation}
where $b= \displaystyle \frac{1}{2} \in \Omega_{1}$.
The operator $A = \displaystyle \frac{\partial^{2}}{\partial x^{2}}$ 
has a complete set of eigenfunctions $(\varphi_{_{k}})$ in $L^{2}(\Omega_{1})$ 
associated with the eigenvalues $(\lambda_{k})$, given by
$$
\varphi_{_{k}}(x) = \sqrt{2} \sin(k \pi x) 
\quad \hbox{and} \quad \lambda_{k} = -k^{2}\pi^{2}
$$
with
$$
S(t)y(x) = \sum_{k=1}^{\infty} e^{\lambda_{k} t} \left\langle y, 
\varphi_{_{k}} \right\rangle_{_{L^{2}(\Omega_{1})}} \varphi_{_{k}}(x).
$$
Then,
$$
H_{\alpha}(t)y(x) = \sum_{k=1}^{\infty} t^{\alpha-1} 
E_{\alpha,\alpha} (\lambda_{k} t^{\alpha})\left\langle y, 
\varphi_{_{k}} \right\rangle_{_{L^{2}(\Omega_{1})}} \varphi_{_{k}}(x),
$$
where $E_{\alpha,\alpha'}(z) := \displaystyle 
\sum_{k=0}^{\infty}\frac{z^{k}}{\Gamma(\alpha k + \alpha')}$, 
\textbf{Re} $\alpha>0$, $\alpha', \, z \in \mathbb{C}$, 
is the generalized Mittag--Leffler function in two parameters (see, e.g., \cite{GO}).
Let $y_{0}(x) = \sin (2\pi x)$ be the initial state to be observed. Then, 
for $\omega_{1} = \displaystyle \left[\frac{1}{6}, \frac{1}{3}\right]$, 
the following result holds.

\begin{Proposition}
There is a state for which the system (\ref{se123})--(\ref{se124}) 
is not weakly observable in $\Omega_{1}$ but it is exactly 
$[\beta_{1}(\cdot), \gamma_{1}(\cdot)]$--observable in $\omega_{1}$.
\end{Proposition}

\begin{proof}
To show that system (\ref{se123})--(\ref{se124}) is not weakly observable in $\Omega_{1}$, 
it is sufficient to verify that $y_{0} \in Ker(K_{\alpha})$.
We have
$$
\begin{array}{rll}
K_{\alpha} y_{0}(x)  
& = & \displaystyle \sum_{i=1}^{\infty} t^{\alpha-1} 
E_{\alpha,\alpha} (\lambda_{i} t^{\alpha})\left\langle y_{0}, 
\varphi_{_{i}} \right\rangle_{_{L^{2}(\Omega_{1})}} \varphi_{_{i}}(b)\\
& = & \displaystyle 2 \sum_{i=1}^{\infty} t^{\alpha-1} 
E_{\alpha,\alpha} (\lambda_{i} t^{\alpha}) \sin\left(\frac{i \pi}{2}\right) 
\int_{0}^{1} \sin(2 \pi x) \sin(i \pi x) dx \\
& = & 0.
\end{array}
$$
Hence, $K_{\alpha} y_{0}(x) = 0$. Consequently, the state $y_{0}$ 
is not weakly observable in $\Omega_{1}$.
On the other hand, one has
$$
\begin{array}{rll}
K_{\alpha} \chi_{_{\omega_{1}}}^{*} \chi_{_{\omega_{1}}} y_{0}(x)  
& = & \displaystyle \sum_{i=1}^{\infty} t^{\alpha-1} 
E_{\alpha,\alpha} (\lambda_{i} t^{\alpha})\left\langle 
\chi_{_{\omega_{1}}}^{*} \chi_{_{\omega_{1}}}y_{0}, 
\varphi_{_{i}} \right\rangle_{_{L^{2}(\Omega_{1})}} \varphi_{_{i}}(b)\\
& = & \displaystyle \sum_{i=1}^{\infty} t^{\alpha-1} 
E_{\alpha,\alpha} (\lambda_{i} t^{\alpha})\left\langle y_{0}, 
\varphi_{_{i}} \right\rangle_{_{L^{2}(\omega_{1})}} \varphi_{_{i}}(b)\\
& = & \displaystyle 2 t^{\alpha-1} E_{\alpha,\alpha} (- \pi^{2} t^{\alpha}) 
\int_{\frac{1}{6}}^{\frac{1}{3}} \sin(2 \pi x) \sin(\pi x) dx \\
& = & \displaystyle \frac{(3\sqrt{3} - 1) t^{\alpha-1} 
E_{\alpha,\alpha} (- \pi^{2} t^{\alpha}) }{6 \pi}\\
& \neq & 0,
\end{array}
$$
which means that the state $y_{0}$ is weakly observable in $\omega_{1}$.
Moreover, for 
$$
\beta_{1} (x) = \left| y_{|\omega_{1}}^{0}(x) \right| - 1 
< y_{|\omega_{1}}^{0}(x) 
$$ 
and 
$$
\gamma_{1} (x) 
= \left| y_{|\omega_{1}}^{0}(x) \right| + 1 
> y_{|\omega_{1}}^{0}(x), \quad \forall x \in \omega_{1},
$$ 
we have $\chi_{_{\omega_{1}}} y_{0}(x) \in [\beta_{1}(\cdot), \gamma_{1}(\cdot)]$ 
and system (\ref{se123})--(\ref{se124}) is exactly 
$[\beta_{1}(\cdot),\gamma_{1}(\cdot)]$-observable in $\omega_{1}$.
The proof is complete.
\end{proof}

Let $\mathcal{G}_{1}$ be the set defined by
$$
\mathcal{G}_{1} = \left\{g \in L^{2}(\Omega_{1}) \; | \; g = 0 \; 
\hbox{in} \; L^{2}(\Omega_{1}) \backslash 
[\beta_{1}(\cdot), \gamma_{1}(\cdot)] \right\}.
$$
From Lemma~\ref{lema5.1}, we see that 
$$
\varphi_{_{0}} \longmapsto \| \varphi_{_{0}} \|_{\mathcal{G}_{1}}^{2} 
= \displaystyle  \int_{0}^{T} \left\| 
C H_{\alpha}(T-t) \varphi_{_{0}} \right\|^{2}  dt
$$
defines a norm on $\mathcal{G}_{1}$. Consider the system
$$
\begin{cases}
\mathcal{Q} \; {}_{t}D_{T}^{\alpha}\Theta(x,t)
=  A^{*} \mathcal{Q}  \Theta(x,t) + \delta(x-b)z(T-t)
& \hbox{in} \quad  \Omega_{1} \times [0,T], \\
\Theta(\xi,t) = 0 & \hbox{on} \quad \partial \Omega_{1} \times [0,T], \\ 
\lim\limits_{t \rightarrow T^{-}}  \mathcal{Q} \; {}_{t}I_{T}^{1-\alpha} \Theta(x,t) 
= 0 &  \hbox{in} \quad \Omega_{1}.
\end{cases}
$$
It follows from Theorem~\ref{thm5.1} that the equation 
$\Lambda : \varphi_{_{0}} \longmapsto \mathcal{P} (\Theta(0))$ 
has a unique solution in $\mathcal{G}_{1}$, which is also 
the initial state $y_{0}$ observed between $\beta_{1}(\cdot)$ 
and $\gamma_{1}(\cdot)$ in the subregion $\omega_{1}$.


\section{Conclusion}
\label{sec:7}

In this paper we have investigated the notion of regional enlarged observability 
for a time fractional diffusion system with Riemann--Liouville fraction derivative 
of order $\alpha \in (0,1]$. We developed an approach that leads to the reconstruction 
of the initial state between two prescribed functions only in an internal subregion 
$\omega$ of the whole domain $\Omega$. We claim that the results here obtained 
can be useful to real problems of engineering. 

As future work, we plan to study problems of regional boundary enlarged observability 
and regional gradient enlarged observability of fractional order distributed parameter systems,
and provide illustrative numerical examples.


\vspace{6pt} 

\authorcontributions{H.Z., A.B. and D.F.M.T. contributed equally to this work.}


\funding{This research was funded by 
\emph{Acad\'emie Hassan II des Sciences et Techniques}, Morocco,
grant number 630/2016 and by
\emph{Funda\c{c}\~ao para a Ci\^encia e a Tecnologia}, Portugal, 
grant number UID/MAT/04106/2013 (CIDMA).}


\acknowledgments{This research was initiated during 
a one-month visit of Zouiten to the Department of Mathematics 
of the University of Aveiro (DMat-UA), November and December 2017, 
followed by a ten days visit to DMat-UA in February 2018.
The hospitality of the host institution and the support 
of Moulay Ismail University, Morocco, and R\&D unit CIDMA, Portugal, 
are here gratefully acknowledged. The authors would like also to thank 
two anonymous reviewers for their critical remarks and precious suggestions, 
which helped them to improve the quality and clarity of the manuscript.}


\conflictsofinterest{The authors declare no conflict of interest.}


\reftitle{References}


\end{document}